\documentclass[letterpaper, 10 pt, conference]{ieeeconf}  

\usepackage[left=0.75in,right=0.75in,top=0.75in,bottom=0.75in]{geometry}

\IEEEoverridecommandlockouts                              

\usepackage{amsmath}
\newcommand*\diff{\mathop{}\!\mathrm{d}}
\newcommand\norm[1]{\left\lVert#1\right\rVert}
\usepackage[ruled,vlined]{algorithm2e}
\usepackage{amssymb}
\usepackage{comment}
\usepackage{color,soul}
\usepackage{mwe}
\usepackage{svg}
\usepackage[hidelinks]{hyperref}
\usepackage{mathtools}
\usepackage{balance}

\usepackage{psfrag}
\usepackage{pstool}
\graphicspath{{Figures/}}

\title{\LARGE \bf Fast and accurate method for computing non-smooth solutions to constrained control problems} 

\author{Lucian Nita$^{1}$, Eduardo M. G. Vila$^{1}$, Marta A. Zagorowska$^{1}$, Eric C. Kerrigan$^{1}$,\\
Yuanbo Nie$^{2}$, Ian McInerney$^{1}$, Paola Falugi$^{1}$

\thanks{*This work has received funding from the Engineering and Physical Sciences Research Council under a Doctoral Training Grant (reference number: EP/T51780X/1), the Active Building Centre project (reference number: EP/V012053/1) and  the Centre for Doctoral Training in High Performance Embedded and Distributed Systems (HiPEDS, Grant Reference EP/L016796/1).}

\thanks{$^{1}$ Department of Electrical \& Electronic Engineering, Imperial College London, SW7~2AZ London, UK, email: 
        {\tt\small \{lucian.nita16, emg216, m.zagorowska, e.kerrigan, i.mcinerney17, p.falugi\}@imperial.ac.uk}}%
\thanks{$^{2}$ Department of Automatic Control and Systems Engineering, University of Sheffield, S1~3JD Sheffield, UK, email:
        {\tt\small y.nie@sheffield.ac.uk}}%
}

\begin{document}

\maketitle
\thispagestyle{empty}
\pagestyle{empty}

\begin{abstract}
Introducing flexibility in the time-discretisation mesh can improve convergence and computational time when solving differential equations numerically, particularly when the solutions are  discontinuous, as commonly found in control problems with constraints.
State-of-the-art methods use fixed mesh schemes, which cannot achieve superlinear convergence in the presence of non-smooth solutions.
In this paper, we propose using a flexible mesh in an integrated residual method. The locations of the mesh nodes are introduced as decision variables, and constraints are added to set upper and lower bounds on the size of the mesh intervals.
We compare our approach to a uniform fixed mesh on a real-world satellite reorientation example. This example demonstrates that the flexible mesh enables the solver to automatically locate the discontinuities in the solution,  has superlinear convergence and faster solve times, while achieving the same accuracy as a fixed mesh.
\end{abstract}

\section{Introduction}
\label{sect:intro}
The differential equations that model a physical system fully describe the dynamical behavior of the system and enable development of control strategies to make the system achieve a desired behavior. Due to the complexity of physical systems, the differential equations describing the dynamics are often impossible to solve analytically, necessitating the use of numerical solvers. The challenge in numerically solving complex differential equations lies in accurately capturing the nature of the solution while using a discrete number of time points, called the time mesh. In this paper, we propose a new method for solving differential equations subject to inequality constraints to a user-specified accuracy by using a flexible time mesh with an integrated residual method.

Differential equations, together with initial values for the states, form an Initial Value Problem (IVP), which can be efficiently solved using explicit~\cite{ahmad_explicit_2004} or implicit methods~\cite{brugnano_blended_2007}. When the differential equations contain boundary conditions, they are called Boundary Value Problems (BVPs), and are customarily solved using collocation methods that compute a polynomial approximation to the solution~\cite{russell_collocation_1972}. Recently, the authors of \cite{nie2022solving} proposed revisiting the use of integrated residual methods for solving dynamic optimization problems, and they highlighted the advantages an integrated residual method has when used to solve complex ordinary differential equations in optimal control. 

The quality of the numerical solution of differential equations depends on the refinement strategy the solver uses to improve the accuracy of the numerical solution. There are three main refinement methods used for improving the solution accuracy~\cite{babuska_h_1992}: (i) \textit{h-methods} that focus on adjusting the size of the time intervals (i.e.\ the time mesh), (ii) \textit{p-methods} that focus on adjusting the polynomial degree, and (iii) \textit{hp-methods} that attempt simultaneous refinement of the time mesh and the degree of the polynomial. State-of-the-art refinement methods are mainly based on adaptive time mesh refinement \cite{liu_adaptive_2015,paiva_adaptive_2015}, but most of the existing approaches are unable to achieve superlinear convergence when discontinuities are present in the solution, and in some cases the numerical solvers may fail to converge to a solution.

The main contribution of this paper is a new method for solving differential equations and control problems. The method is based on an integrated residual method with a flexible time mesh, which allows one to accurately capture potential discontinuities in the solution. Our approach introduces the time-mesh nodes as variables in an optimisation problem while numerically integrating the residual function over each time interval. Using a least squares approach to compute the residuals enables us to solve a wide range of problems including ordinary differential equations (ODEs) and differential algebraic equations (DAEs), BVPs, and consistently over-determined systems. This new method can also be used to handle the dynamic constraints in optimal control problems (OCPs), which often have discontinuous solutions. Additionally, the proposed method can handle differential inclusions or discontinuities that may appear in the dynamics.
The proposed flexible time mesh provides better performance than a fixed time mesh, with a numerical example showing superlinear convergence and an improved solution accuracy, while also decreasing the computational time.

In Section~\ref{sect:pbdef}, we present an overview of the integrated residual method for solving differential equations and extend the method to the solution of control problems. We introduce our flexible time mesh scheme in Section~\ref{sect:qflex}. In Section~\ref{sect:results}, we use our proposed flexible mesh scheme to solve a minimum-time satellite reorientation problem. Section~\ref{sect:concl} concludes the paper and presents an outline for future works.

\section{Problem Definition}
\label{sect:pbdef}
\subsection{Differential Equations and Control Problems}
Many control problems can be formulated as finding the solution to one or more feasibility problems of the form
\begin{subequations}
\label{eq:Feasibility}
\begin{align}
    \text{find} \quad & x(\cdot), u(\cdot)\\
    \text{s.t.} \quad & F(\dot{x}(t), x(t), u(t), t) = 0, \quad \forall t \in \mathcal{T}\label{eq:mainEq}\\
    & G(\dot{x}(t), x(t), u(t), t) \leq 0, \quad \forall t \in \mathcal{T},
\end{align}
where $\mathcal{T} := [t_0, t_f] \subset \mathbb{R}$ is the time interval over which the problem is defined, $x : \mathbb{R} \rightarrow \mathbb{R}^{N_x}$ are the state differential variables and are forced to be continuous, $\dot{x} : \mathbb{R} \rightarrow \mathbb{R}^{N_x}$  are the time derivatives of the state, $u : \mathbb{R} \rightarrow \mathbb{R}^{N_u}$ are the algebraic variables, which can model the control inputs, and are allowed to be discontinuous. The function $F : \mathbb{R}^{N_x} \times \mathbb{R}^{N_x} \times \mathbb{R}^{N_u} \times \mathbb{R} \rightarrow \mathbb{R}^{N_F}$, which is typically called the dynamics function, defines a set of $N_F$ equality constraints that the controlled system have to satisfy and $G : \mathbb{R}^{N_x} \times \mathbb{R}^{N_x} \times \mathbb{R}^{N_u} \times \mathbb{R} \rightarrow \mathbb{R}^{N_g}$ defines $N_g$ path inequality constraints. Moreover, the problem may include boundary constraints
\begin{align}
    \Psi_E(x(t_0), x(t_f), t_0, t_f) & = 0, \label{eq:fixedpt}\\
    \Psi_I(x(t_0), x(t_f), t_0, t_f) & \leq 0, 
    \label{eq:fixedpt2}
\end{align}
where 
$\Psi_E : \mathbb{R}^{N_x} \times \mathbb{R}^{N_x} \times \mathbb{R} \times \mathbb{R} \rightarrow \mathbb{R}^{N_E}$ are the boundary equality constraints, and $\Psi_I : \mathbb{R}^{N_x} \times \mathbb{R}^{N_x} \times \mathbb{R}  \times \mathbb{R} \rightarrow \mathbb{R}^{N_I}$ are the boundary inequality constraints. 
\label{eq:dop}
\end{subequations} 

Differential equations are usually included in problem~\eqref{eq:dop} as  equality constraints in~\eqref{eq:mainEq}. 
Note that solutions to~\eqref{eq:Feasibility} are \emph{non-unique} and \emph{non-smooth}, in general, even if $F$, $G$ $\Psi_E$ and $\Psi_I$ are all smooth functions.

Of particular interest is to note that 
problems of the form~\eqref{eq:dop} need to be solved in certain classes of direct methods for solving optimal control problems~\cite{nie2022solving}. 
\label{sec:de}
\subsection{Discretisation and Numerical Solution}
The problem introduced in Section \ref{sec:de} is an infinite-dimensional problem because of the dependence on time. To numerically solve an infinite-dimensional feasibility problem using finite-dimensional optimisation methods, the problem is discretized. The goal of the resulting problem is to find approximating functions $\tilde{x} : \mathbb{R} \rightarrow \mathbb{R}^{N_x}$
and $\tilde{u} : \mathbb{R} \rightarrow \mathbb{R}^{N_u}$, such that the integrated square of the residual 
\begin{equation}
    \epsilon_R:=\frac{1}{(t_f-t_0)\cdot N_F}\int_{t_0}^{t_f}{\norm{ F(\dot{\tilde{x}}(t),\tilde{x}(t),\tilde{u}(t),t)}^2_2}\diff t
    \label{eq:approxErr}
\end{equation}
is minimised. The residual $\|F(\dot{\tilde{x}}(t),\tilde{x}(t),\tilde{u}(t),t)\|^2_2$ captures how close the numerical solution $(\tilde{x}(t),\tilde{u}(t))$ is to an exact solution $(x(t),u(t))$ at each time instant $t$, since~\eqref{eq:mainEq} becomes zero at the exact solution. A scaling factor is added in order to average the residual over the domain~$\mathcal{T}$ and over all components of the dynamics equations. A necessary and sufficient condition for achieving convergence of $(\tilde{x}(\cdot),\tilde{u}(\cdot))$  to an exact solution $(x(\cdot),u(\cdot))$ is to make the integrated residual go to zero. Since the exact solution $(x(\cdot),u(\cdot))$ is non-unique, in general, defining an error metric based on the exact solution (such as $\| \tilde{x}(t)-x(t) \|^2_2$ for example) would not be practical.

The approximating functions $\tilde{x}$ and $\tilde{u}$ are typically represented using polynomial basis functions \cite[Sect.~1.17.1]{betts_practical_2010}. Finding these approximating functions consists of determining a finite number of coefficients for the polynomials. In order to numerically enforce constraints in \eqref{eq:mainEq}, we insert the approximating functions $\tilde{x}, \tilde{u}$ into function $F$. In the least squares approach we are proposing here, $F(\cdot)$ is replaced with a function of the coefficients of the polynomials. A procedure to obtain the approximating function and approximate solutions to the  differential equations will be discussed in Section~\ref{sec:IRM}.

The choice of the basis functions for $(\tilde{x}(\cdot),\tilde{u}(\cdot))$ determines the possible solution space, which implies that the exact solution $(x(\cdot),u(\cdot))$ may not be representable in that solution space. This means the integrated residual $\epsilon_R$ may be non-zero in general, and the only case when $\epsilon_R$ can be zero is if the span of the
chosen basis functions contains an element of the solution
set. When basis functions span the solution set, it is then possible to numerically determine the exact  trajectories.

\begin{figure}[tb]
\centering
\includegraphics[width=0.48\textwidth]{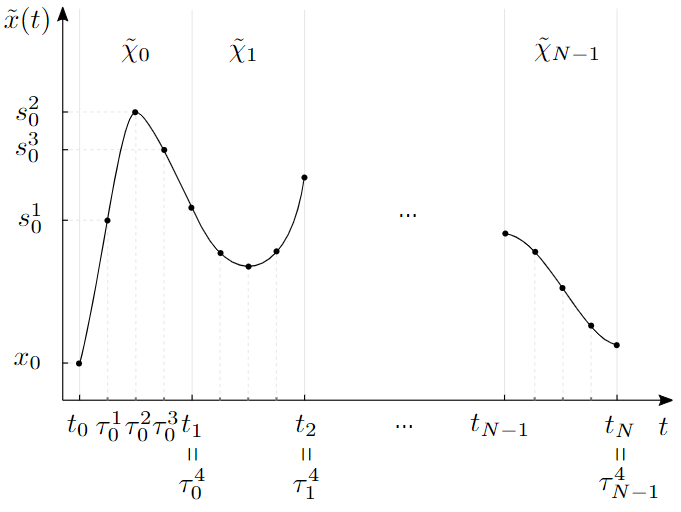}
\caption[Notation]{Approximation method and notation used in constructing a numerical solution to problem \eqref{eq:dop}.}
\label{fig:NotationArrows}
\end{figure}

To reduce the approximation error in the numerical computation of \eqref{eq:approxErr} while maintaining numerical stability, the interval $[t_0, t_f]$ is split into $N$ subintervals 
\begin{equation}
    \label{eq:tidef}
    [t_{i}, t_{i+1}] =: \mathcal{T}_i\subset\mathcal{T} \quad \forall i \in \{0, 1, \dots,N-1\}
\end{equation}
 delimited by the time-mesh points $t_i$ with $t_{N} = t_{f}$. A diagram showing the subintervals and the basis functions can be seen in Figure~\ref{fig:NotationArrows}, with each subinterval $\mathcal{T}_i$ having a component of $\tilde{x}$ given by the basis function $\tilde{\chi}_{i}$ with 5 support points $\tau_{i}^{j}, j \in \{0, 1, 2, 3, 4\}$.


In our implementation, we use a piecewise polynomial of degree $a$ chosen from the Lagrange basis 
\begin{equation}
\begin{aligned}
    \tilde{x}(t) & =\tilde{\chi}_i(t) \quad \forall t\in[t_i,t_{i+1}]\\
    \tilde{\chi}_i(t) & \coloneqq \sum_{j=0}^{a}{s_i^j \cdot \prod_{\substack{k=0 \\ k\neq j}}^{a}{\frac{t-\tau_i^k}{\tau_i^j-\tau_i^k}}},
\end{aligned}
\label{eq:approxfct}
\end{equation}
where $\tau_i^j\in [t_i,t_{i+1}], \forall j \in \{0,\dots,a\}$ are the supports for the interpolating polynomial in the $i$-th interval, and $s_i^j$ for $j\in\{0,\dots,a\}$ are the decision variables that need to be obtained from the optimisation process for each interval $i\in\{0,\dots,N-1\}$. Since the basis function is an interpolating polynomial, the $s_i^j$ are the $\tilde{\chi}_i$ function values when evaluated at a finite number of time points $\tau_i^j$, giving
\begin{equation}
    s_i^j\coloneqq\tilde{\chi}_i(\tau_i^j)=\tilde{x}(\tau_i^j), \ \forall i \in \{0,\dots,N-1\}, 
    j \in \{0,\dots,a\}. 
\end{equation}

To enforce state continuity at mesh nodes, we introduce the constraints
\begin{subequations}
\label{eq:constra}
\begin{align}
    \tilde{\chi}_i(t_{i+1})&=\tilde{\chi}_{i+1}(t_{i+1})\quad \forall i\in \{0,\dots,N-2 \},\\
    \tilde{\chi}_0(t_0)&=x_0,\\
    \tilde{\chi}_{N-1}(t_f)&=x_f.
\end{align}
\end{subequations}

One way of ensuring the constraints in \eqref{eq:constra} are satisfied is by placing internal supports at the interval boundaries $\tau_i^0=t_i$ and $\tau_i^a=t_{i+1} \quad \forall i \in \{0,\dots,N-1\}$. With the internal supports at the interval boundaries, the same decision variables $s_i^a=s_{i+1}^0$ can be used to represent both $\tilde{\chi}_i(t_{i+1})$ and $\tilde{\chi}_{i+1}(t_{i+1}), \forall i \in \{0,\dots,N-2\}$.

A similar strategy is used for transforming a control input~$u(\cdot)$ into a piecewise polynomial $\tilde{u}(\cdot)$ based on a finite number of decision variables $c_i^j$ where $i$ is the index of the subinterval and $j \in \{0, 1, ..., b\}$ is the $j^\text{th}$ polynomial coefficient. Note that the degree of polynomial approximation may be different for the states and inputs ($a\neq b$), so the supports may not overlap. In this case, we will evaluate all the functions at the supports used by the highest degree polynomial approximation. For the functions which do not have a support present at the point of evaluation, we will use the interpolated value. Finally, $\dot{\tilde{x}}(\cdot)$ is obtained by differentiating the state approximation function with respect to time. 

A typical design choice is to fix the values $t_i$ in the time-mesh before solving problem \eqref{eq:dop} to create a fixed mesh. The most common approach is to place mesh nodes $t_i$ at equidistant locations to create a uniform mesh. However, there are also other options for selecting the predefined position of the nodes \cite{betts_mesh_1998}. A uniform time-mesh is given by
\begin{equation}
    t_i = t_0 + i\cdot\frac{t_f-t_0}{N},
    \label{eq:unimesh}
\end{equation} 
where $t_{i}$ are the time-mesh nodes and $i \in \{0,\dots,N\}$ is the time-mesh node index.

The choice of the number of intervals for a fixed mesh depends on the required solution accuracy with respect to the error metric $\epsilon_R$ from \eqref{eq:approxErr}. An in-depth discussion on the influence of $N$ on the error metric can be found in Section~\ref{sect:qflex}.

\label{sec:numsol}
\subsection{Integrated Residual Methods}
\label{sec:IRM}

Integrated residual methods (IRMs) can be used to solve differential equations of the form described in \eqref{eq:mainEq}. Additionally, they can be used to solve  optimal control problems with inequality constraints using the method in~\cite{nie2022solving}. Classical collocation methods enforce constraint \eqref{eq:mainEq} only at a finite number of points, called collocation points. A  drawback of the collocation approach is that one does not have direct control over the numerical error encountered in-between the collocation points. In contrast, IRMs aim to enforce constraint~\eqref{eq:mainEq}  by directly minimising the integrated residual over the entire desired time-span $\mathcal{T}$.     

To efficiently solve feasibility and control problems with discontinuous 
solutions, which are otherwise difficult to tackle, we will use an integrated residual method for solving differential equations that is similar to \cite{eason_review_1976} and \cite{mortari_high_2019}. The idea of this approach is to transform the feasibility problem \eqref{eq:dop} into a constrained minimisation problem of the form 
\begin{subequations}
    \label{eq:leastsq}
\begin{align}
    \min_{\tilde{x}(\cdot),\tilde{u}(\cdot)} \quad & \epsilon_{R}\label{eq:objlsq}\\ 
    \textrm{s.t.} \quad & G(\dot{\tilde{x}}(t), \tilde{x}(t), \tilde{u}(t), t) \leq 0, && \forall t \in \mathbf{\tau} \label{eq:varplsq}\\
    & \Psi_E(\tilde{x}(t_0), \tilde{x}(t_f), t_0, t_f) = 0, \label{eq:fixedplsq}\\
    & \Psi_I(\tilde{x}(t_0), \tilde{x}(t_f), t_0, t_f) \leq 0, \label{eq:infixedplsq}
    \end{align}
\end{subequations}
with the integrated residual~\eqref{eq:approxErr} as the objective to be minimised. To keep the implementation simple, path inequality constraints are only enforced at the support time points \begin{align}
    \mathbf{\tau}:=\{\tau_i^j\mid i=0,1,\dots N-1; j=0,1,\dots a\}.  
\end{align}

In practice, the integrals in \eqref{eq:approxErr} and \eqref{eq:objlsq} are solved numerically using quadrature rules. In our implementation we used Gaussian quadrature with $Q$ quadrature points \cite{venkateshan_chapter_2014}. The integral in \eqref{eq:approxErr} is computed numerically and since $F$ is a general nonlinear function, this will introduce an additional numerical error, namely the \emph{quadrature error}

\begin{equation}
    \epsilon_Q\coloneqq\Big|\epsilon_R-\sum_{i=0}^{N-1}{\sum_{k=1}^{Q}{w_i^k \cdot \norm{ F(\dot{\tilde{x}}(\rho_i^{k}),\tilde{x}(\rho_i^{k}),\tilde{u}(\rho_i^{k}),\rho_i^{k})}^2_2\Big|}}
    \label{eq:quaderr}
\end{equation} where $w_i^k$ for $k\in\{1,\dots,Q\}$ are the $Q$ Gaussian quadrature weights associated with the integration interval $[t_i,t_{i+1}]$ and~$\rho_i^{k}$ are the quadrature nodes for the interval $[t_i,t_{i+1}]$. After introducing the quadrature rule, the discretisation method from Section~\ref{sec:numsol} can be used to turn problem \eqref{eq:leastsq} into the discrete formulation
\begin{equation}
    \begin{aligned}
    \min_{\mathbf{s},\mathbf{c}} \quad & \sum_{i=0}^{N-1}{\sum_{k=1}^{Q}{w_i^k \cdot \norm{ F(\dot{\tilde{x}}(\rho_i^{k}),\tilde{x}(\rho_i^{k}),\tilde{u}(\rho_i^{k}),\rho_i^{k})}^2_2}}\\
        \textrm{s.t.} \quad & \eqref{eq:constra}, \eqref{eq:varplsq}, \eqref{eq:fixedplsq}, \eqref{eq:infixedplsq},
    \end{aligned}
    \label{eq:numsollsq}
\end{equation} 
where $\mathbf{s} \in \mathbb{R}^{N(a+1)}$ and $\mathbf{c}\in \mathbb{R}^{N(b+1)}$
are vectors consisting of all the  polynomial coefficients  $s_i^j, c_i^j$. 
Problem \eqref{eq:numsollsq} can now be solved directly using available NLP solvers. 

To make sure the quadrature error is negligible, we need to check that the number of quadrature points $Q$ is sufficiently high. This can be done after the original solve by recomputing $\epsilon_R$ (to ensure the numerical integration has converged) using a  larger $Q$ than was used to solve~\eqref{eq:numsollsq} (e.g.\ with $2Q$ quadrature points) and evaluating $\epsilon_Q$ according to \eqref{eq:quaderr}. If the difference between the converged $\epsilon_R$ and the objective value obtained from optimization problem \eqref{eq:numsollsq} is above a certain tolerance $\varepsilon_{\text{quad,tol}}$ (i.e.\ $\epsilon_{Q} > \varepsilon_{\text{quad,tol}}$), this implies that the quadrature error is significant and the solution needs to be recomputed using more quadrature points $Q$.

\section{Improving Accuracy Through Mesh Design Updates}
\label{sect:qflex}
\subsection{Flexible mesh}
\label{sect:fm}
When selecting a suitable number of mesh intervals $N$, one should recall that the solution $(x(\cdot),u(\cdot))$ is in a particular space, which may be different from the one spanned by the selected basis functions of the discretisation. Thus, the residual for the discretisation might be non-zero.

In order to improve the solution accuracy, the conventional mesh refinement process relies on increasing the number of time nodes or changing the polynomial degree according to a predefined algorithm based on error metrics~\cite{betts_practical_2010}. One popular strategy is to start with a coarse mesh (small $N$), evaluate the solution accuracy with respect to the error metric~$\epsilon_R$ while ensuring the quadrature error $\epsilon_Q$ remains below a threshold $\varepsilon_{\text{quad,tol}}$, and then recompute a new solution for a finer mesh (large $N$) if the error is above a certain tolerance level.

For many classes of problems, mesh refinement can be ineffective and inefficient. Consider for example a problem with a discontinuous solution with discontinuities at unknown locations. Unless it happens by chance for a node to be placed exactly at a discontinuity, we can consider without loss of generality that the discontinuity lies inside the interval $(t_{i},t_{i+1})$. This means that one is trying to approximate a discontinuous function $x(\cdot)$ by a continuous function $\tilde{\chi}_i(\cdot)$. Increasing the polynomial degree will not help, since a ringing phenomenon will start to occur. Equivalently, increasing the number of mesh nodes will not improve the solution for the interior of the interval where the discontinuity is located. 

To automate this process of designing the time mesh and capturing discontinuities and regions of high gradients more accurately, as well as to improve the overall solution for a given number $N$ of mesh intervals, we propose to include the mesh nodes $t_i$ as decision variables in our optimization problem. 


To ensure the quadrature error $\epsilon_Q$ remains insignificant as the mesh nodes move, we introduce a rule restricting the allowable change in the individual mesh interval length. Consider a flexibility parameter $\phi \in [0,1)$, then the constraints 
\begin{subequations}
    \begin{align}
        t_{i+1}-t_i & \leq (1+\phi)\cdot\frac{t_f-t_0}{N},\   \forall i\in\{1,\dots,N-1\}
    \label{eq:quasi_meshu}\\
    t_{i+1}-t_i & \geq (1-\phi)\cdot\frac{t_f-t_0}{N},\   \forall i\in\{1,\dots,N-1\}
            \label{eq:quasi_meshl}
    \end{align}
        \label{eq:quasi_mes}%
\end{subequations}
will be included in the optimisation problem. It follows that if $\phi\in(0,1)$ then $t_0<t_1<\dots<t_{N-1}<t_N$ and hence no overlaps are possible.

Hence, the optimal node locations can be determined automatically by solving the following optimisation problem
\begin{equation}
    \begin{aligned}
    \min_{\mathbf{s},\mathbf{c},\mathbf{t}} \quad & \sum_{i=0}^{N-1}{\sum_{k=1}^{Q}{w_i^k \cdot \norm{ F(\dot{\tilde{x}}(\rho_i^{k}),\tilde{x}(\rho_i^{k}),\tilde{u}(\rho_i^{k}),\rho_i^k)}^2_2}}\\
    \textrm{s.t.} \quad & \eqref{eq:constra}, \eqref{eq:varplsq}, \eqref{eq:fixedplsq}, \eqref{eq:infixedplsq},
    \eqref{eq:quasi_mes}
    \end{aligned}
    \label{eq:qfexmesh}
\end{equation}
where $\mathbf{t} \in \mathbb{R}^{N-1}$ is the vector of mesh points $t_i \quad\forall i \in \{1,\dots,N-1\}$ where the initial and final times ($t_0$ and $t_f$) are excluded, since they are fixed. Note that quadrature points $\rho_i^k$ and internal supports $\tau_i^j$ become a function of $t_i$ and $t_{i+1}$, hence
 they are shifted and scaled versions of the original values. 

\subsection{Error and Performance Considerations}

When implementing the flexible mesh together with the integrated residual approach, we ensure the accuracy for the trajectories $\tilde{x}(\cdot)$ and $\tilde{u}(\cdot)$ in-between mesh nodes. Compared to classical time-marching schemes (shooting methods) or point-wise residual minimisation (collocation), integrated residual methods have the benefit of producing a solution with a more uniform error over the whole time domain. 



When comparing the problem formulations in \eqref{eq:numsollsq} and~\eqref{eq:qfexmesh}, note that the only difference is represented by the variables we are optimising over. 
Observe also that $\phi=0$ leads to the fixed uniform mesh solution, i.e.\ $t_{i+1}-t_i = (t_f-t_0)/N$. 
As a result, problem~\eqref{eq:qfexmesh} includes the fixed uniform mesh problem \eqref{eq:numsollsq}.
By enlarging the solution space, successfully solving \eqref{eq:qfexmesh} with a flexible mesh should provide a solution with a residual error that is not worse than the one obtained using a fixed mesh in \eqref{eq:numsollsq}, provided the quadrature error is sufficiently small.

For a given value for the flexibility parameter $\phi$, an algorithm that implements~\eqref{eq:qfexmesh} simplifies the mesh refinement procedure, since the mesh node locations are automatically determined and the main parameters that can be changed are the number of nodes $N$ and the polynomial degree. This  increases the accuracy of the obtained solution $\tilde{x}(\cdot)$ and $\tilde{u}(\cdot)$. Additionally, the proposed strategy speeds up the computation of a solution, as demonstrated in Section~\ref{sect:results}.

\section{Results}
\label{sect:results}
An integrated residual solver was developed in the Julia~v1.6 programming language, implementing fixed and flexible mesh schemes for the solution of constrained DAEs. Dynamic variables were parameterised by Lagrange polynomials in the barycentric form~\cite{berrut_barycentric_2004}. These were discretised across~$N$ time intervals, as illustrated in Figure~\ref{fig:NotationArrows}. Residuals were integrated with Gauss-Legendre quadrature of a sufficiently high order. Derivatives were evaluated using forward and reverse automatic differentiation. The least squares method described in Section~\ref{sec:IRM} was applied using Ipopt~\cite{wachter_implementation_2006} as the NLP solver (relative convergence tolerance set to $10^{-8}$). All tests were performed on an Intel\textsuperscript{\textregistered} Core\textsuperscript{\texttrademark} i7-1065G7 at 1.3\,GHz with 16\,GB of RAM.
\subsection{Satellite Reorientation Example}
The NASA X-ray Timing Explorer spacecraft is modelled as a 3D rigid body with moments of inertia $I_{xx} > I_{yy} > I_{zz}$. The orientation is defined by quaternions $q = [q_1, q_2, q_3, q_4]$, where $\norm{q}_2 = 1$ must always hold. The spacecraft must be reoriented ($\Delta q_1 = 150^{\circ}$) in minimum time, subject to the saturation of control torques $\norm{u}_{\infty} \leq 50$\,N\,m. To model the dynamics, the least squares method allows for the DAE formulation to be used to directly enforce the quaternion magnitude, as opposed to the ODE formulation, which requires high-accuracy integration, hence we used the DAE formulation~\cite[Sect.~6.8]{betts_practical_2010}. 

\begin{figure}
    \centering
    \includegraphics[width=\columnwidth]{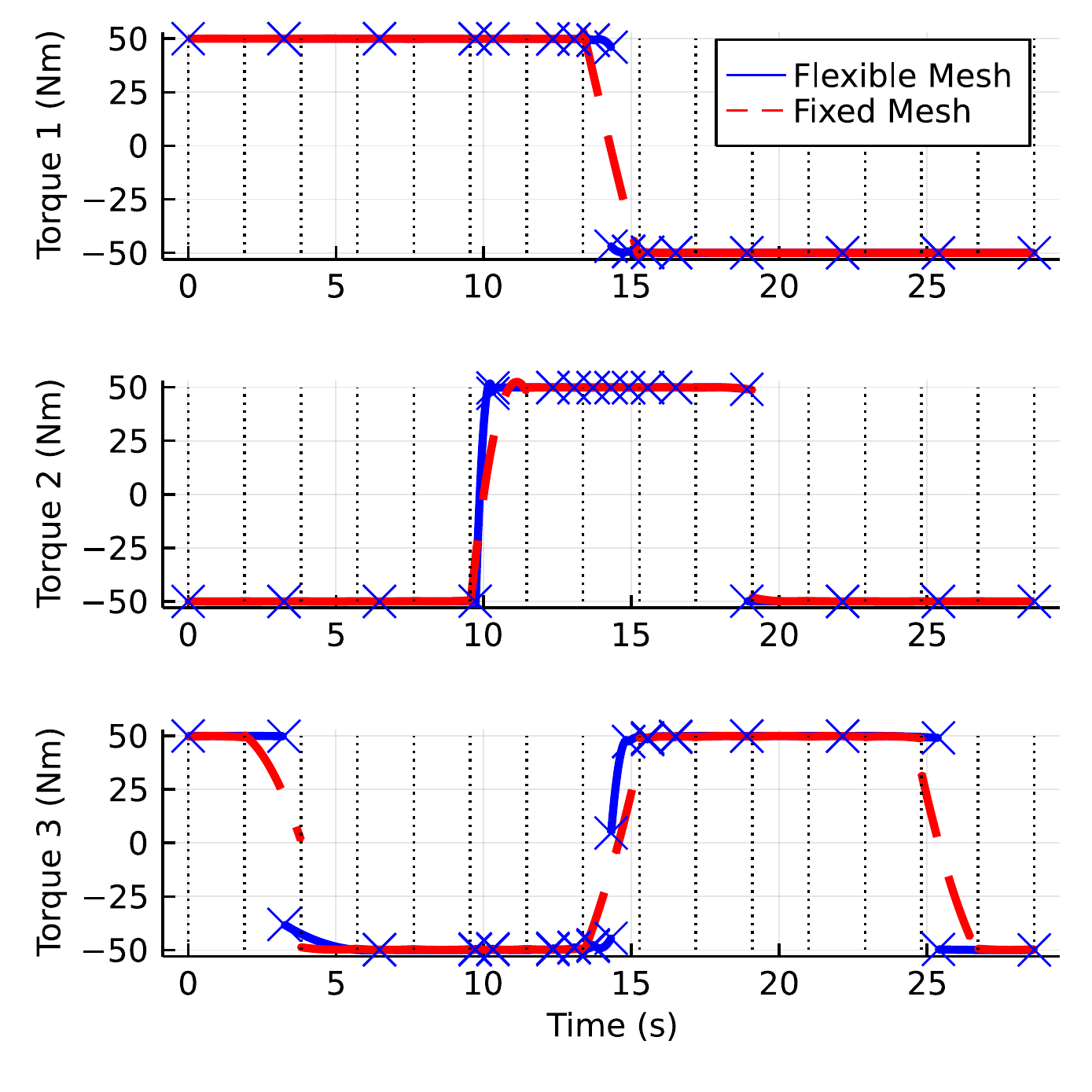}
    \caption{Control solution to satellite reorientation using fixed and flexible meshes with 15 intervals. Dotted vertical lines indicate the location of the fixed mesh points, blue crosses indicate the location of the flexible mesh points.}
    \label{fig:solution}
\end{figure}

Despite the fact that the dynamic functions are continuous (albeit nonlinear), this  control problem yields a discontinuous solution with switches at non-trivial times. Figure~\ref{fig:solution} shows a solution to~\eqref{eq:dop} when the final time is set to the optimum $t_f = 28.630408\,\mathrm{s}$~\cite[Sect.~6.8]{betts_practical_2010}.
To perform the maneuver in minimum time, control torques 2 and 3 ($u_2$ and $u_3$, respectively) introduce coupled rotations, reducing the effective moment of inertia.

As expected, the flexible mesh is able to capture the discontinuities by shifting the mesh points to the switch times. On the other hand, the fixed mesh has to compromise accuracy at discontinuities, with the (accidental) exception of control torque 1 ($u_1$) where the uniform distribution \eqref{eq:unimesh} has a time node coinciding with the discontinuity at $t = t_f/2$.
\subsection{Convergence}
\begin{figure}
    \centering
    \includegraphics[width=\columnwidth]{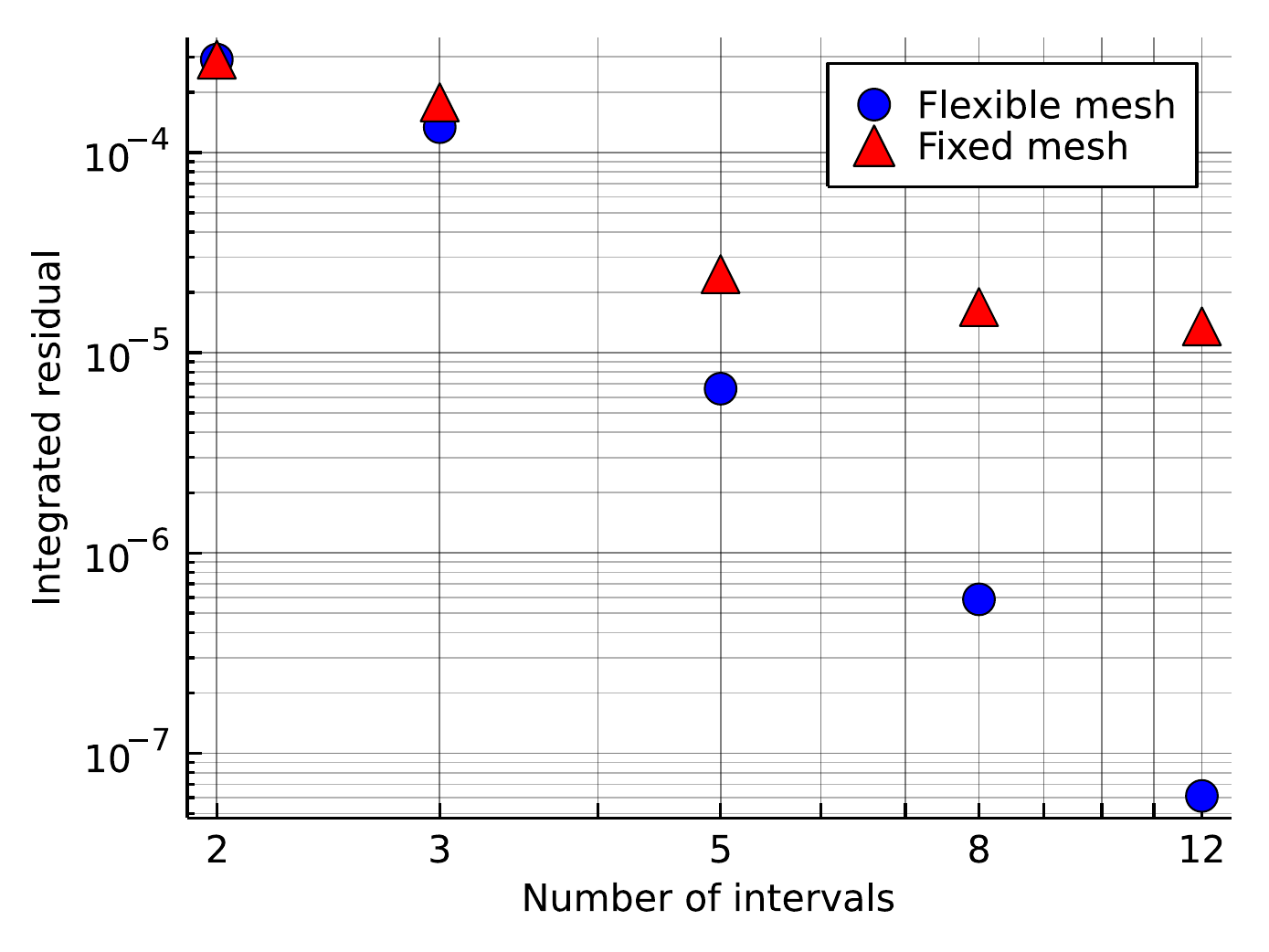}
    \caption{Log-log plot of integrated residual for solutions to the satellite reorientation problem for both mesh schemes.}
    \label{fig:convergence}
\end{figure}

To assess the order of convergence for the proposed scheme, the problem was solved for an increasing number of intervals using both flexible and fixed mesh points. The initial guess for the NLP time-mesh decision variables was set to be the uniform mesh, while the state and control decision variables $\mathbf{s}, \mathbf{c}$ were cold started (matrices full of zeros were given as the initial guess). Figure~\ref{fig:convergence} shows the value of the integrated residuals plotted against $N$. 
The slopes of linear fits for $N \geq 5$ report orders of convergence of 0.70 and 5.55 for fixed and flexible meshes, respectively. Lagrange polynomials of degree 4 were used as basis functions, with a sufficiently high quadrature order, $Q = 7$.

Once the number of intervals is larger than the number of discontinuities, the flexible scheme is able to locate the position of the discontinuities. Hence, for a larger number of intervals, the flexible mesh is able to achieve a higher order of convergence than the fixed mesh. 
\subsection{Computation time}
\begin{figure}
    \centering
    \includegraphics[width=\columnwidth]{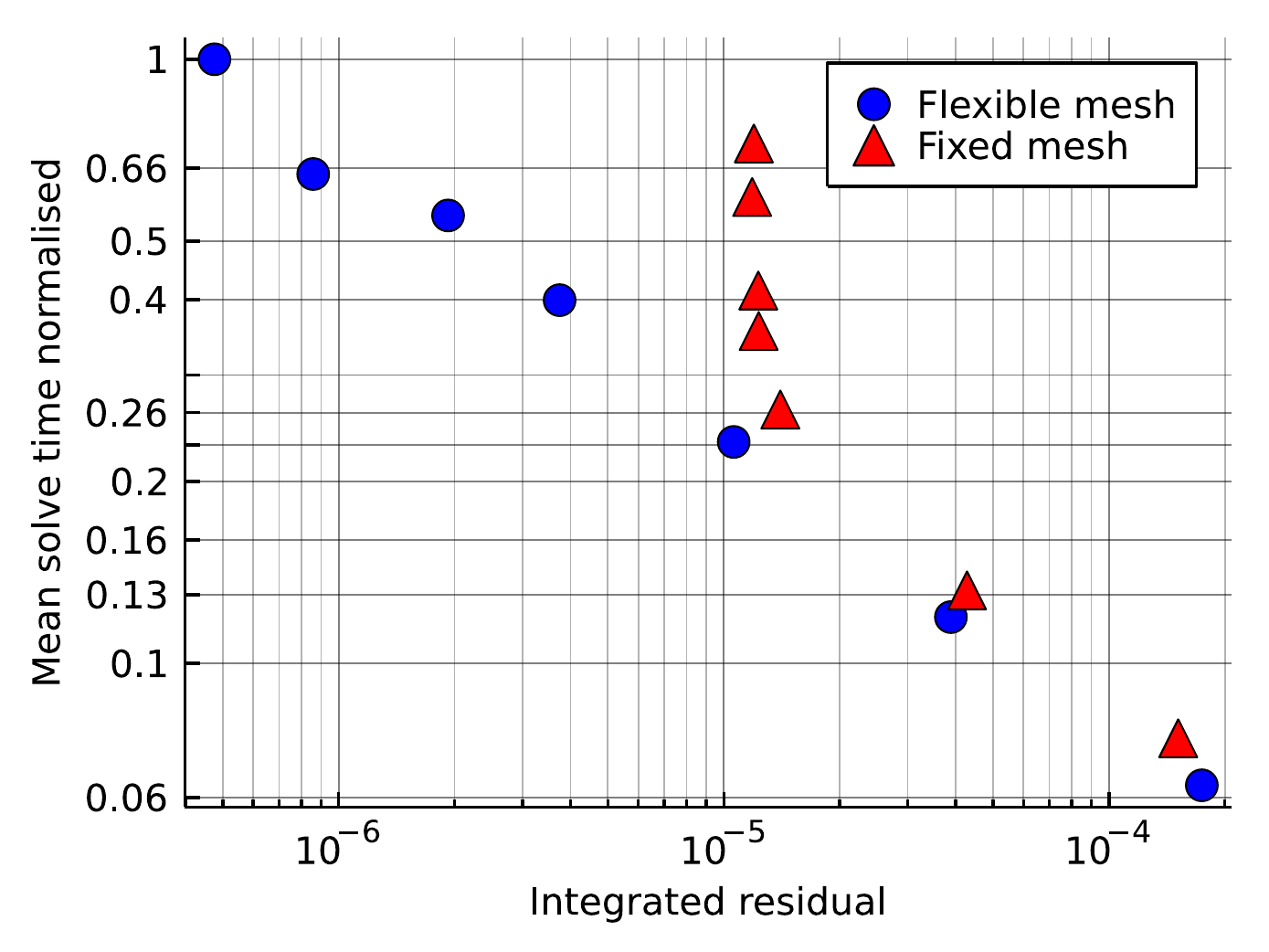}
    \caption{Log-log plot of performance against integrated residuals for both mesh schemes with different number of intervals $N$ between 2 and 24 and $\phi=0.5$. The solve times are normalised by the highest value.}
    \label{fig:performance}
\end{figure}

Adding flexible mesh nodes increases the size of the resulting NLP, and introduces nonlinearities. On the other hand, this reduces the number of NLP iterations to achieve a given integrated residual. 
Computation time was assessed by comparing the solution time averaged with 5 repeated solves. Figure~\ref{fig:performance} shows that for low accuracy solutions ($\epsilon_R > 10^{-5}$), there is little difference in computational time between the two schemes. However, the flexible mesh is able to achieve a significantly better accuracy given the same solve time as a fixed mesh. Conversely, faster solve times can be achieved for a given accuracy.






\section{Conclusions and future works}
\label{sect:concl}
Finding discontinuous solutions to constrained control problems is a challenge, especially if a fixed time-mesh is used. In  this  paper, we  proposed using  a  flexible  mesh  in  an  integrated  residual  method to capture possible discontinuities in the solution. By including the mesh points as decision variables in the integrated residual method, we are able to capture discontinuities in the solution, while preserving accuracy of the solution. The proposed scheme was tested on a challenging aerospace satellite reorientation problem with nonlinear dynamics and a discontinuous solution. In contrast to using a fixed mesh, the new method with a flexible mesh successfully identified the location of the  discontinuities, provided enough mesh intervals were chosen. The flexible scheme showed super-linear convergence once all discontinuities are captured, whereas the fixed mesh scheme plateaued.
Compared to a fixed mesh scheme, the underlying NLP solver converged to a solution in fewer iterations, resulting in the faster overall computation of a solution at a higher accuracy, despite the added overhead of augmenting the underlying optimisation problem.

These preliminary numerical results on order of convergence and computational performance are encouraging. Future work could include developing formal proofs on the order of convergence. The use of a flexible mesh should also be compared numerically against state-of-the-art adaptive mesh refinement methods, which do not include mesh points as decision variables, but iteratively add fixed mesh points to intervals with a large residual. Similar sophisticated mesh refinement schemes could also be explored when using a flexible mesh.

The flexible mesh method presented in this paper can be adapted to solve optimal control problems using an iterative algorithm. This can be done, for example, by solving a sequence of feasibility problems~\eqref{eq:qfexmesh} with an increasing number of mesh intervals $N$ until~$\epsilon_R$ is below a given threshold. The optimal control problem is then solved by defining and solving another suitable optimisation problem  as 
in~\cite{nie2022solving}. 


\balance

\bibliographystyle{plain}
\bibliography{root.bib}

\end{document}